\documentclass[11pt]{article}
\usepackage[utf8]{inputenc}
\usepackage{tikz}
\usepackage{graphicx}
\usepackage{hyperref,verbatim}
\usepackage{enumitem,bm}
\usepackage{todonotes}
\usepackage[font=footnotesize,labelfont=bf]{caption}
\usepackage{placeins}
\usepackage{float}
\usepackage{setspace}

\usepackage{pgf,tikz,subcaption}
\usepackage{mathtools}
\usetikzlibrary{arrows,shapes}
\usetikzlibrary{decorations.pathreplacing}
\usetikzlibrary{calc,fit}
\usepackage{fullpage}

\usepackage{amsthm,amsmath,amssymb}
\usepackage{url}
\usepackage[capitalise]{cleveref}

\newtheorem{conj}{Conjecture}
\newtheorem{corollary}[conj]{Corollary}
\newtheorem{theorem}[conj]{Theorem}
\newtheorem{lem}[conj]{Lemma}

\newtheorem{prop}[conj]{Proposition}

\theoremstyle{definition}

\newtheorem{example}[conj]{Example}

\newtheorem{remark}[conj]{Remark}

\newcommand*{\myproofname}{Proof}

\newcommand{\C}{\mathcal{C}}

\tikzstyle{uStyle}=[shape = circle, minimum size = 6.5pt, inner sep = 0pt,
outer sep = 0pt, draw, fill=white, semithick]
\tikzstyle{sStyle}=[shape = rectangle, minimum size = 4.5pt, inner sep = 0pt,
outer sep = 0pt, draw, fill=white, semithick]
\tikzstyle{lStyle}=[shape = circle, minimum size = 4.5pt, inner sep = 0pt,
outer sep = 0pt, draw=none, fill=none]
\tikzset{every node/.style=uStyle}

\definecolor{MainPart}{RGB}{221,235,247}
\definecolor{MainBorder}{RGB}{93,130,166}
\definecolor{CoreCol}{RGB}{117,107,177}
\definecolor{DummyCol}{RGB}{166,189,219}
\definecolor{AuxPart}{RGB}{224,243,219}
\definecolor{AuxBorder}{RGB}{49,130,116}
\definecolor{XCol}{RGB}{102,194,165}
\definecolor{YCol}{RGB}{255,237,111}
\definecolor{ZCol}{RGB}{244,109,67}

\tikzset{
  vertex/.style={circle, draw=black, fill=white, line width=.25pt, inner sep=0pt, minimum size=4.8pt},
  corev/.style={vertex, fill=CoreCol!75},
  dummyv/.style={vertex, fill=DummyCol!95},
  auxv/.style={vertex, fill=white},
  mainpart/.style={rounded corners=3pt, draw=MainBorder, fill=MainPart, line width=.75pt, minimum width=2.75cm, minimum height=4.25cm},
  corebox/.style={rounded corners=2pt, draw=CoreCol!65!black, fill=CoreCol!13, line width=.45pt, minimum width=.78cm, minimum height=1.34cm},
  dummyset/.style={rounded corners=2pt, draw=DummyCol!72!black, fill=DummyCol!50, line width=.48pt, minimum width=.54cm, minimum height=.98cm},
  auxpart/.style={rounded corners=3pt, draw=AuxBorder, fill=AuxPart, line width=.75pt, minimum width=2.75cm, minimum height=3.20cm},
  xset/.style={rounded corners=2pt, draw=XCol!65!black, fill=XCol!45, line width=.48pt, minimum width=1.64cm, minimum height=.52cm},
  yset/.style={rounded corners=2pt, draw=YCol!55!black, fill=YCol!62, line width=.48pt, minimum width=1.64cm, minimum height=.46cm},
  zset/.style={rounded corners=2pt, draw=ZCol!70!black, fill=ZCol!40, line width=.48pt, minimum width=1.64cm, minimum height=.52cm},
  pathEdge/.style={draw=black!62, line width=.45pt},
  complete/.style={draw=black!78, line width=2.05pt, line cap=round, opacity=.84},
  legendbox/.style={rounded corners=1.5pt, draw=black!25, fill=white, inner sep=3pt}
}

\newcommand{\MainPart}[2]{%
  \begin{scope}[shift={(#2,0)}]
    \node[mainpart] (M#1) at (0,0) {};
    \node[font=\scriptsize\bfseries, text=MainBorder!85!black] at (0,1.92) {$M_{#1}$};

    \node[corebox] (CoreBox#1) at (0,1.08) {};
    \node[corev] (c#1a) at (0,1.42) {};
    \node[corev] (c#1b) at (0,1.08) {};
    \node[corev] (c#1c) at (0,0.74) {};

    \foreach \r/\x in {1/-0.82,2/0,3/0.82}{
      \node[dummyset] (D#1\r) at (\x,-0.95) {};
      \node[dummyv] at ($(D#1\r.center)+(0, .24)$) {};
      \node[dummyv] at ($(D#1\r.center)+(0,-.24)$) {};
    }
  \end{scope}%
}

\newcommand{\AuxPartNode}[4]{% name, x, y, label
  \begin{scope}[shift={(#2,#3)}]
    \node[auxpart] (A#1) at (0,0) {};

    \node[xset] (X#1) at (0,.86) {};
    \foreach \xx in {-.56,-.28,0,.28,.56}{\node[auxv] at (\xx,.86) {};}

    \node[yset] (Y#1) at (0,.02) {};
    \node[auxv] at (0,.02) {};

    \node[zset] (Z#1) at (0,-.86) {};
    \foreach \xx in {-.40,0, .40}{\node[auxv] at (\xx,-.86) {};}

    \node[font=\scriptsize\bfseries, text=AuxBorder!85!black] at (0,-1.34) {#4};
  \end{scope}%
}

\newcommand{\affiliation}{\footnote}
\makeatletter
\def\@fnsymbol#1{\ensuremath{\ifcase#1\or *\or \dagger\or \ddagger\or \mathsection\or \|\or **\or \dagger\dagger \or \ddagger\ddagger \else\@ctrerr\fi}}
\makeatother

\title{Supersaturation and reconfiguration diameter of independent transversals}
\date{\today}

\author{
Wouter Cames van Batenburg\affiliation{D\'epartement d'Informatique, Universit\'e libre de Bruxelles, Belgium, \textsf{\href{mailto:w.p.s.camesvanbatenburg@gmail.com}{w.p.s.camesvanbatenburg@gmail.com}}.}
\and
Penny Haxell\affiliation{Department of Combinatorics and Optimization, University of Waterloo,
Canada, \textsf{\href{mailto:pehaxell@uwaterloo.ca}{pehaxell@uwaterloo.ca}}.}
\and
Ross J. Kang\affiliation{Korteweg--de Vries Institute for Mathematics, University of Amsterdam,
Netherlands, \textsf{\href{mailto:r.kang@uva.nl}{r.kang@uva.nl}}.}
\and
Ronen Wdowinski \affiliation{Institute of Discrete Mathematics, Graz University of Technology,
Austria, \textsf{\href{mailto:wdowinski@math.tugraz.at}{wdowinski@math.tugraz.at}}.}
}

\begin{document}
\maketitle

\begin{abstract}
Consider the collection of independent transversals of a vertex-partitioned graph of maximum degree $\Delta$.
It is a classic result that when all parts have size at least $2\Delta$, the collection is nonempty.
Recently, it was found that if the parts have size at least $2\Delta+1$, then every pair of independent transversals is 
% or: "then every two independent transversals are"
connected, within the collection, by a sequence of one-vertex reconfigurations.
We prove, among other results, that there are exponentially many independent transversals when the parts have size at least $2\Delta$, and that the reconfiguration diameter is at most linear in the number of parts when the parts have size at least $4\Delta$.
\end{abstract}

\section{Introduction}

Given an integer $\Delta\ge 3$, let $G$ be a graph of maximum degree at most $\Delta$. Consider an arbitrary vertex-partition $\mathcal P=\{P_1,\ldots,P_n\}$ for $G$.
An \emph{independent transversal} of $(G,\mathcal{P})$ is an independent set of $G$ containing exactly one vertex from each part.
This is guaranteed to exist, by a result of Haxell~\cite{haxell1995condition}, if all parts of $\mathcal P$ have size at least $2\Delta$. The size-constraint $2\Delta$ is precisely optimal, in that there are examples of $(G,\mathcal{P})$ where the parts are only of size
$2\Delta-1$ for which there is no independent transversal~\cite{szabo2006Transversals} (see also~\cite{haxell2024constructing}). 

Consider the collection of all independent transversals of $(G,\mathcal{P})$. %, especially under variation of the size-constraint for the parts of $\mathcal P$.
How does this collection change as we increase the size-constraint beyond $2\Delta$?
\begin{enumerate}
\item\label{q:supersaturation}
Can we guarantee it has a certain cardinality, as a function of natural parameters in $(G,\mathcal{P})$?
\item\label{q:geometry}
What about changes in the collection's qualitative structure, such as its topology or geometry?
\end{enumerate}
The result of~\cite{haxell1995condition} guarantees only that the cardinality in Question~\ref{q:supersaturation} is at least $1$ for any size-constraint at least $2\Delta$, and so Question~\ref{q:supersaturation} constitutes a natural `supersaturation' refinement.
To make Question~\ref{q:geometry} precise, 
%we augment the collection by a graph with which we define the distance between two independent transversals of $(G,\mathcal{P})$.
%The \emph{IT-reconfiguration graph} of $(G,\mathcal{P})$ is the graph whose vertices are the independent transversals of $(G,\mathcal{P})$ where two that differ on exactly one part are joined by an edge. (We may often drop the `IT-' prefix where it causes no confusion.)
%The reconfiguration problem asks when two independent transversals can be transformed into one another by changing one chosen vertex at a time, while remaining an independent transversal after each step.
%Beyond existence, we are interested in what the entire space of independent transversals looks like. How ``close'' are two independent transversals to each other? To make that question precise, 
we study  
the \emph{(independent transversal) reconfiguration graph} of $(G,\mathcal{P})$, which is the graph whose vertices are the independent transversals of $(G,\mathcal{P})$ and whose edges join two transversals that differ in exactly one part. Under what circumstances is the reconfiguration graph connected, and how small is its diameter? These kinds of questions are relevant for applications to sampling and counting algorithms via Markov chains (see, e.g., \cite{jerrum1995very}).

Buys, Kang, and Ozeki~\cite{buys2025reconfiguration} used combinatorial methods to prove that the reconfiguration graph of $(G,\mathcal{P})$ is connected when every part has size at least $2\Delta+1$. Wdowinski reproved this result within a framework of topological methods~\cite{wdowinski2025hall}, analogous to those that were used to prove Haxell's existence theorem (e.g., \cite{aharoni2006intersection, aharoni2007independent, meshulam2003domination}). If every part has size exactly $2\Delta$, connectivity can fail only if some parts $\mathcal{P}' \subseteq \mathcal{P}$ induce $|\mathcal{P}'|$ disjoint copies of $K_{\Delta,\Delta}$ in $G$.
Wdowinski~\cite{wdowinski2026tight}, resolving a question in~\cite{buys2025reconfiguration}, exactly characterized which partitions may occur in the exceptional instances.
These results give a sharp qualitative picture for connectivity around Haxell's threshold, while leaving quantitative diameter and supersaturation questions open. 

%More recently, Cambie, Cames van Batenburg, Cranston, van den Heuvel, and Kang proved an $O(|V(G)|^2)$ transformation result for list-colourings under the local condition $|L(v)|\ge \deg(v)+1$ in the non-frozen setting, and showed that reducing the list size of even one vertex to $\deg(v)$ can shatter the reconfiguration space~\cite{cambie2025reconfiguration}.

\subsection{Diameters of reconfiguration graphs}

One of our main results provides a linear upper bound on the diameter under the stronger assumption that the parts have size $|P_i|\ge 4\Delta$.
The diameter bound $2n-1$ is sharp due to examples coming from correspondence-colouring reconfiguration~\cite[Proposition~2(ii)]{cambie2024optimally}. 

\begin{theorem}\label{thm:linear-diameter}
Let $G$ be a graph of maximum degree $\Delta$, and let $\mathcal{P}=\{P_1,\dots,P_n\}
$ be a partition of $V(G)$ such that
$|P_i|\ge 4\Delta$
for every $i\in[n]$.
Then the reconfiguration graph of independent transversals of $(G,\mathcal{P})$ has diameter at most
$2n-1$.
\end{theorem}

Related diameter questions have been studied before for list colourings, which can be viewed as independent transversals of an auxiliary vertex-partitioned graph, and hence constitute a special case (see \cite{haxell2001note}). In that setting, a list colouring is already guaranteed to exist when parts have size $\Delta+1$, and the associated reconfiguration graph is connected when the parts have size at least $\Delta+2$. At that same threshold of parts of size $\Delta+2$, Jerrum~\cite{jerrum1995very} gave an argument that yields a $2^{O(n)}$ upper bound on the diameter of the reconfiguration graph, Cereceda~\cite{cereceda2007diss} improved it to $O(n^2)$, and finally Cambie, Cames van Batenburg, and Cranston~\cite{cambie2024optimally} showed that the diameter is at most $2n-1$. See~\cite{CvB26, DeMeyer26} for a generalization to hypergraphs and refinements of this sharp bound.\\
%and conjectured that it is at most $\lfloor{3n/2\rfloor}$, which would be optimal.
%for all values of $n$.
%For reconfiguration of \emph{correspondence colourings} (which form a generalization of list colourings and a special case of independent transversals) they showed the upper bound $2n-1$ holds as well and is sharp for all values of $n$.\\ 
%2n-1 is  only optimal for correspondence colouring; for list colouring it is conjecturally at most n + matching number.

For general independent transversals, this raises the question of whether the $2n-1$ bound in Theorem~\ref{thm:linear-diameter} might still hold close to the $2\Delta$ connectivity threshold. We show that this is the case asymptotically if the average degree of the vertices in each block is $o(\Delta)$ (see Corollary \ref{cor:average-degree-linear}), but this still leaves open the general case. Is the diameter still at most $2n-1$, or at least linear in $n$, under the weaker condition that all parts have size at least $2\Delta+1$? And when the parts have size $\ge 2\Delta$, does every connected component have diameter $O(n)$? So far, the best result we know \cite{hajebi2026algorithms} is that when parts have size at least $2\Delta+1$, the diameter is polynomial in $n$ for every fixed $\Delta$, but the dependence on $\Delta$ in the exponent is huge. In this current work, we show that linear diameter fails below the connectivity threshold, and that connected components can have at least quadratic diameter. 

\begin{theorem}\label{thm:simulationC3Ps}
%For every $\Delta\ge 3$ and infinitely many $n$, there is a graph $G$ of maximum degree at most $\Delta$, 
%The previous above formulation could be misinterpreted as the weaker statement that the "infinitely many n" purely arise from the infinitely many possible values of \Delta.
For every $\Delta\ge 3$ there are infinitely many $n$ for which there exists a graph $G$ of maximum degree at most $\Delta$, 
together with a partition $\mathcal P$ of $V(G)$ into $n$ parts each of size $2\Delta-1$, 
such that the reconfiguration graph of independent transversals of $(G,\mathcal{P})$ has a connected component of diameter at least
$
\frac{n^2}{4\Delta^2}.
$
\end{theorem}

Our construction simulates the $3$-colouring reconfiguration graph of a path (from \cite{bonamy2011diameter}) inside the independent transversal reconfiguration graph.

\subsection{Number of independent transversals}

In the direction of counting independent transversals, our next main result is that there exist exponentially many independent transversals when the parts have size $|P_i| \ge 2\Delta$, where the exponential dependence is in the number of parts $n$. 

\begin{theorem} \label{thm:exponentially-many}
Let $G$ be a graph of maximum degree $\Delta$, and let $\mathcal{P}=\{P_1,\dots,P_n\}$ be a partition of $V(G)$ such that $|P_i|\ge 2\Delta$. Then the number of independent transversals of $(G,\mathcal{P})$ is at least $(2\Delta)^{\left(1-\frac{1}{2\Delta-1}\right)\frac{n}{2}}$.
\end{theorem}

Our proof is an interesting application of Lindsey's edge-isoperimetric inequality \cite{lindsey1964assignment} on the Hamming graph. Our theorem appears to be the first result demonstrating exponentially many independent transversals at the threshold $2\Delta$ that works for all values of $n$ and $\Delta$. 

Exponential lower bounds on the number of independent transversals have previously been obtained through the Lov\'asz Local Lemma and variations of it, but they require larger part sizes. For example, Wanless and Wood~\cite{wanless2022general} proved the existence of at least $(t/2)^n$ independent transversals when $|P_i| = t \ge 4\Delta$ for all $i$. This was actually shown with the maximum degree $\Delta$ replaced by the \emph{maximum block average degree} $D$ of $(G,\mathcal{P})$, which is defined as the maximum over all parts $P_i$ of the average degree of the vertices in $P_i$. In this average degree setting, the part size threshold $4D$ was shown to be best possible even for the existence of a single independent transversal \cite{groenland2023graphs,haxell2024degree}. And while it is suggestive based on Theorem \ref{thm:linear-diameter} that the reconfiguration graph of $(G,\mathcal{P})$ ought to be connected when $|P_i| \ge 4D$, that hunch fails spectacularly: It is not hard to construct examples with two parts $P_1,P_2$ of \emph{arbitrarily large} size and $D<2$ such that the independent transversal reconfiguration graph is disconnected.
%Indeed. Let the two parts be $P_1=\{x_1,\ldots,x_t\}$ and $P_2=\{y_1,\ldots, y_t\}$.Put edges $x_1y_j$ and $x_iy_1$ for all $i,j \ge 1$.Then the number of edges is $2t_2$ so the maximum block average degree is $D=(2t-2)/2 < 2$, but the independent transversal $\{x_1,y_1\}$ forms an isolated vertex in the reconfiguration graph. 

Another exponential-type lower bound is due to Tang and Zhao \cite{tang2025number}. Their main result implies that if $|P_i| \ge 2\Delta$ for all $i$, then the number of independent transversals is at least $c_n \Delta^n$, where $c_n > 0$ depends only on $n$. However, they require that $\Delta$ be sufficiently large compared to $n$. We are more interested in the situation when $\Delta$ is fixed and $n$ is growing. Buys, Kang, and Ozeki \cite{buys2025reconfiguration} posed the following conjecture giving a sharp lower bound on the number of independent transversals at Haxell's theshold $2\Delta$. This would be tight by a simple construction (see Example \ref{ex:cyclic}).

\begin{conj} \label{conj:supersaturation}
Let $G$ be a graph with maximum degree $\Delta$, and let $\mathcal{P}=\{P_1,\ldots,P_n\}$ be a partition of $V(G)$ with $|P_i| = 2\Delta$ for all $i$. Then the number of independent transversals of $(G, \mathcal{P})$ is at least $2 \Delta^n$.
\end{conj}

For comparison, the total number of (not necessarily independent) transversals of $\mathcal{P}$ is $(2\Delta)^n$. Our last main result is proving the first non-trivial case of Conjecture \ref{conj:supersaturation}, namely when $n = 3$ (Theorem \ref{thm:three-parts}). Interestingly, while there are various tight constructions for Conjecture \ref{conj:supersaturation} when $n \in \{1,2\}$, we show that the mentioned tight construction (Example \ref{ex:cyclic}) is unique when $n = 3$, and we conjecture that this is the case for all $n \ge 3$.

\subsection{Organization}
This paper is organized as follows. In Section \ref{sec:linear-diameter}, we prove Theorem \ref{thm:linear-diameter} on linear diameter. In Section \ref{sec:quadratic-diameter}, we prove Theorem \ref{thm:simulationC3Ps} on quadratic diameter components. In Section \ref{sec:exponentially-many}, we prove Theorem \ref{thm:exponentially-many} as well as a more elementary result coming from double-counting. In Section \ref{sec:small-number-of-parts}, we prove Conjecture \ref{conj:supersaturation} when $n \le 3$.

\section{Linear diameter when parts have size $\ge 4\Delta$} \label{sec:linear-diameter}

In this section, we prove Theorem \ref{thm:linear-diameter} demonstrating that the reconfiguration graph of $(G,\mathcal{P})$ has linear diameter when the parts of $\mathcal{P}$ have size at least $4\Delta$. 
Before that, let us recall Haxell's Theorem \cite{haxell1995condition,haxell2001note}.

\begin{theorem}[Haxell~\cite{haxell1995condition, haxell2001note}] \label{thm:existence}
Let $G$ be a graph of maximum degree at most $\Delta$ with a vertex partition $\mathcal{P}=\{P_1,\ldots, P_n\}$. If $|P_i| \ge 2\Delta$ for every $i$,
then $(G,\mathcal{P})$ has an independent transversal. In fact, every vertex is contained in some independent transversal.
\end{theorem}

Theorem \ref{thm:linear-diameter} will be a direct corollary of Theorem~\ref{thm:existence} and the following theorem, which more generally provides a linear upper bound on the diameter under the condition that the existence of an independent transversal is sufficiently robust under vertex removal.

\begin{theorem}\label{thm:main_extended}
Let $G$ be a graph with maximum degree at most $\Delta$ and let $\mathcal{P}=\{P_1,\dots,P_n\}$ be a partition of $V(G)$  into $n$ parts. Suppose that whenever one removes at most $2\Delta$ vertices from each part of $\mathcal{P}$, the resulting induced subgraph $G'$ and vertex partition $\mathcal{P}'$ has an independent transversal. Then the reconfiguration graph of independent transversals of $(G,\mathcal{P})$ has diameter at most
$2n-1$.
%Suppose that for every induced subgraph $G'$ with vertex partition $\mathcal{P}'$ that is obtained by removing at most $2\Delta$ vertices from each part of $\mathcal{P}$, $(G',\mathcal{P}')$ has an independent transversal.  Then the reconfiguration graph of independent transversals of $(G,\mathcal{P})$ has diameter at most $2n-1$.
\end{theorem}

%The proof of Theorem~\ref{thm:main_extended} uses a suitable independent transversal $C$ such that both $A$ and $B$ can be reconfigured to $C$ in at most $n$ steps.
\begin{proof}
Let $A=\{a_1,\ldots, a_n\}$ and $B=\{b_1,\ldots, b_n\}$ be two independent transversals of $(G,\mathcal{P})$. We will show that $
\operatorname{dist}(A,B)\le 2n-1.$ Essentially, the proof proceeds by finding some suitable independent transversal $C$ such that both $A$ and $B$ can be reconfigured to $C$ in at most $n$ steps.
We construct an ordering $\rho$ of the parts $P_i$, together with large subsets
$R_i\subseteq P_i$, with the ordering $\rho$ chosen so that every vertex of $R_i$ avoids all vertices $a_j$ and $b_j$ appearing
later than $i$ in the order $\rho$. We construct $\rho$ greedily. 

Suppose that $U\subseteq [n]$ is the set of indices of parts that are not yet ordered. For
$i\in U$, define
$$
\operatorname{Block}_i(U):=
P_i\cap N\bigl(\{a_j \cup b_j :j\in U\setminus\{i\}\}\bigr).
$$
Thus $\operatorname{Block}_i(U)$ is the set of vertices in $P_i$ blocked by the still-unordered vertices of $A$ and $B$. Since
$\sum_{i\in U}|\operatorname{Block}_i(U)|
\le
\sum_{j\in U}\deg(a_j) + \deg(b_j)
\le
|U|\cdot 2\Delta,
$
some $i\in U$ satisfies
$$
|\operatorname{Block}_i(U)|\le 2\Delta.
$$
Put this part $P_i$ next in the order $\rho$, and define
$
R_i:=P_i\setminus \operatorname{Block}_i(U)$.
Then
$$
|R_i|\ge |P_i|-2\Delta.
$$
It follows from the hypothesis of Theorem~\ref{thm:main_extended}
that the graph induced by
$\bigcup_{i=1}^n R_i$
has an independent transversal
$C=\{c_i:c_i\in R_i\}
$.
Moreover, by construction, every vertex of $R_i$ is nonadjacent to every $a_j$ and $b_j$ whose part appears after $P_i$ in
the order $\rho$. \\

We now show that $A$ can be reconfigured to $C$ in at most $n$ steps. We reconfigure the parts in the order
$\rho$. Suppose we are reconfiguring the part $P_i$. At this moment, the parts earlier than $P_i$ in $\rho$
have already been changed to their vertices of $C$, while the parts later than $P_i$ in $\rho$ are still
at their vertices of $A$.
We now reconfigure $a_i$ to $c_i$. This is allowed because:

\begin{itemize}
    \item $c_i$ is nonadjacent to all earlier vertices of $C$, since $C$ is independent; and
    \item $c_i$ is nonadjacent to all later vertices $a_j$ of $A$, since $c_i\in R_i$.
\end{itemize}

By the same argument, the transversal $B$ can be reconfigured
to $C$ in at most $n$ steps according to the same order $\rho$.
This already gives that
$\operatorname{dist}(A,B)\le \operatorname{dist}(A,C)+ \operatorname{dist}(B,C)\le n+n=2n.$
However, we can do slightly better.
Let $r \in \{1,\ldots, n\}$ be the index whose part $P_r$ comes last in the order $\rho$.
As above we reconfigure $a_i$ and $b_i$ to $c_i$ according to order $\rho$, except we do not yet execute the last reconfiguration, from $a_r$ to $c_r$ and from $b_r$ to $c_r$. Instead we can then finish by directly reconfiguring $a_r$ to $b_r$, saving one reconfiguration step.
\end{proof}

As another application of Theorem \ref{thm:main_extended}, we prove Corollary \ref{cor:average-degree-linear} below, showing that we have linear diameter whenever the parts have size $(2+o(1))\Delta$, under the additional assumption that the maximum block average degree $D$ of $(G,\mathcal{P})$ is $o(\Delta)$. Recall that the \textit{maximum block average degree} of $(G,\mathcal{P})$ is defined as $D \coloneqq \max_{P_i \in \mathcal{P}} \frac{1}{|P_i|} \sum_{v \in P_i} d_G(v)$. In place of Haxell's Theorem \ref{thm:existence}, we use the following theorem that comes from variations of the Lov\'asz Local Lemma \cite{kang2022colorings, wanless2022general} (see also \cite{bissacot2011improvement}). 

\begin{theorem}[Kang, Kelly \cite{kang2022colorings}; Wanless, Wood \cite{wanless2022general}] \label{thm:existence-average-degree}
Let $G$ be a graph and let $\mathcal{P}=\{P_1,\ldots, P_n\}$ be a vertex partition. If $(G,\mathcal{P})$ has maximum block average degree $D$ and $|P_i| \ge 4D$ for every $i$, then $(G,\mathcal{P})$ has an independent transversal.
\end{theorem}

\begin{corollary} \label{cor:average-degree-linear}
Let $G$ be a graph of maximum degree $\Delta$ and maximum block average degree $D$, together with a vertex partition
$\mathcal{P}=\{P_1,\ldots,P_n\}$ such that $|P_i|\ge 2\Delta+6\sqrt{D\Delta}$
 for every $i\in [n]$.
Then the reconfiguration graph of independent transversals of
\((G,\mathcal P)\) has diameter at most \(2n-1\).
%Same proof works with $|P_i|\ge 2\Delta+ (2+2\sqrt{3}) \sqrt{D\Delta}$ instead.
\end{corollary}
\begin{proof}
We verify the hypothesis of Theorem~\ref{thm:main_extended}.
%We may assume $D\ne 0$ because otherwise the result is immediate.
Let $G'$ be obtained from $G$ by deleting
at most $2\Delta$ vertices from each part, and put $P_i'\coloneqq P_i\cap V(G')$ for every $i$.
Then $|P_i'|\ge 6\sqrt{D\Delta}$.
Let $D'$ be the maximum block average degree of $(G',\{P_1',\ldots,P_n'\})$, say that the average degree is maximised at $P_j'$. Since $D\le \Delta$, we have
 $$D'=\frac{1}{|P_j'|}\sum_{v\in P_j'} d_{G'}(v)
 \le
 \frac{1}{|P_j'|}\sum_{v\in P_j} d_G(v)
 \le
 \frac{|P_j|D}{|P_j|-2\Delta}
 \le
 \frac{(2\Delta+6\sqrt{D\Delta})D}{6\sqrt{D\Delta}}
 \le
 \frac{4}{3}\sqrt{D\Delta}.$$
 Hence $4D' <6\sqrt{D\Delta}\le |P_i'|$
for every \(i\). By Theorem~\ref{thm:existence-average-degree}, it follows that $(G',\{P_1',\ldots, P_n'\})$ has an independent transversal. Therefore, Theorem~\ref{thm:main_extended} applies to give the diameter result.
\end{proof}

\begin{remark}
Due to the construction in \cite[Proposition~2(ii)]{cambie2024optimally}, for every $\Delta=n-1$ and $p \ge 3$ there exist examples of graphs $G$ with maximum degree $\Delta$ and a partition $\mathcal{P}=\{P_1,\ldots, P_n\}$ where $|P_i|\ge p$ for all $i$ and the reconfiguration graph of $(G,\mathcal{P})$ has diameter at least $2n-1$. Thus the upper bound $2n-1$ in Theorem~\ref{thm:linear-diameter} is sharp for all values of $\Delta$.
\end{remark}

%\begin{remark}
%With weaker sufficient conditions for the existence of an independent transversal, 
% the condition $|P_i|\geq 4\Delta$ in Theorem~\ref{thm:main} might be lowered to something like $|P_i| \geq (3+o(1)) \Delta$. Not clear yet how to lower the condition closer to $|P_i|\ge 2\Delta+1$.
%\end{remark}

\begin{remark}
If the parts have size $k\ge 2\Delta+1$, then there exists at least one ``unfrozen'' part. More precisely:
Suppose $\mathcal{P}=\{P_1,\ldots,P_n\}$ and $|P_i|=k$ for all $i$.
Given independent transversals $A$ and $B$, say that $P_i$ is \emph{$(A,B)$-unfrozen} if some $y\in P_i$ satisfies
$y\not\sim A\setminus\{a_i\}$ and $y\not\sim B\setminus\{b_i\}$.
Suppose for a contradiction that no such $y$ exists for any $i$. Then with
$e_A(i)=e(A\setminus\{a_i\},P_i)$ and $e_B(i)=e(B\setminus\{b_i\},P_i)$,
we have $e_A(i)+e_B(i)\ge k$ for all $i$.
Since $\sum_i e_A(i)\le n\Delta$ and $\sum_i e_B(i)\le n\Delta$, it follows that $kn \le 2n\Delta$, contradiction. So there must exist some $(A,B)$-unfrozen part $P_i$. This could be a somewhat desirable property for a potential proof by induction, because it means that $a_i$ can be reconfigured to $y$ without being blocked by either $A$ or $B$.\end{remark}

\section{Quadratic diameter components when parts have size $\le 2\Delta-1$} \label{sec:quadratic-diameter}

In this section, we prove Theorem~\ref{thm:simulationC3Ps}. When the parts of a partition have size less than $2\Delta$, there might not be an independent transversal, and if there does exist one, the reconfiguration graph may not be connected. However, here we zoom in on the situation that the reconfiguration graph is nonempty, and we ask how large the maximum diameter among all connected components can be. We show that it can be at least $\frac{n^2}{4\Delta^2}$. We use the following result as input. 

\begin{theorem}[Bonamy, Johnson, Lignos, Patel, Paulusma~\cite{bonamy2011diameter}]\label{thm:pathreconfiguration}
For every $s\ge 2$, the $3$-colouring reconfiguration graph of the $s$-vertex path $P_s$ has diameter at least $\frac{s^2-1}{4}$.
\end{theorem}

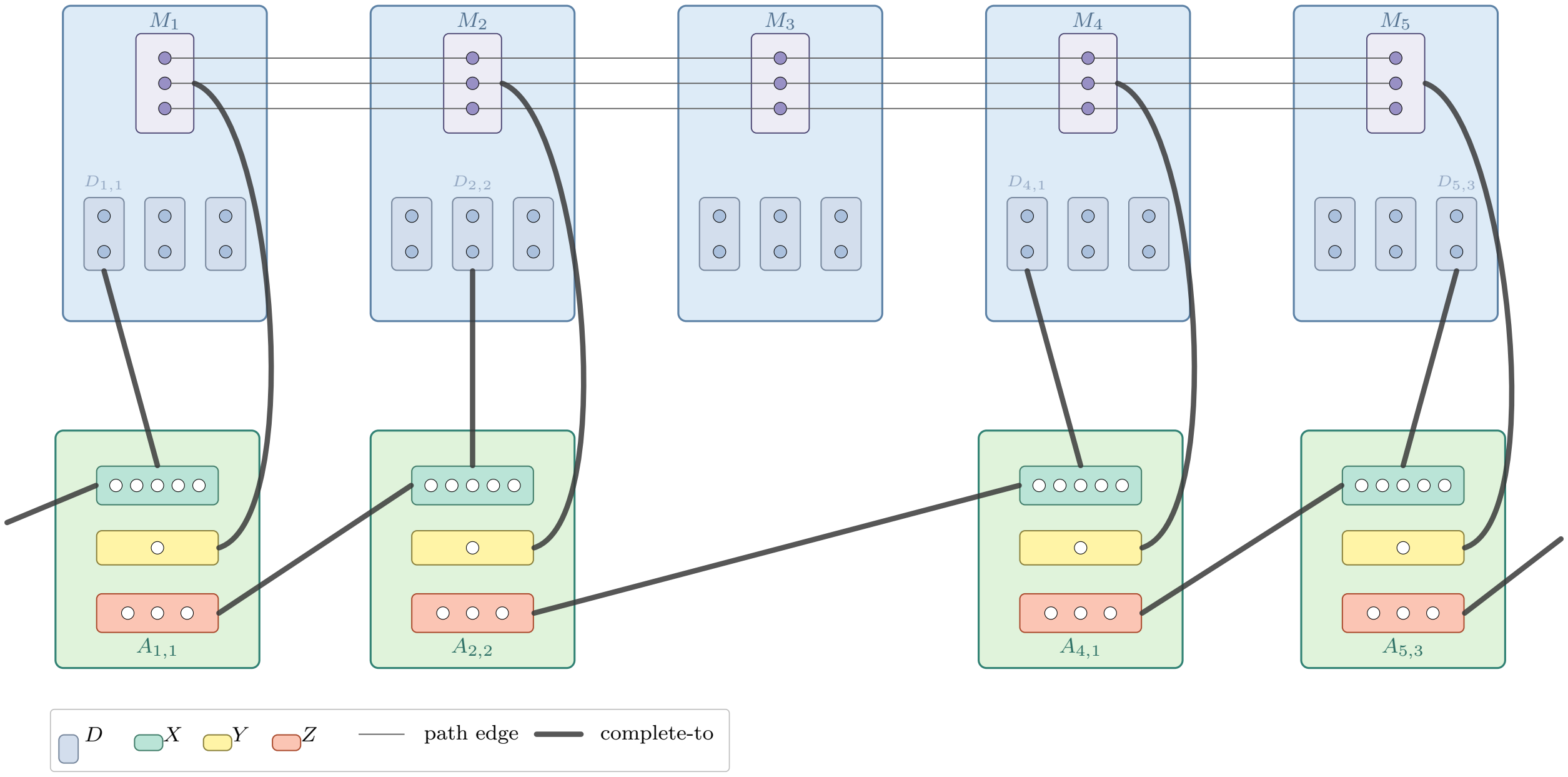
\begin{figure}[H]
\centering
\resizebox{\textwidth}{!}{%
\begin{tikzpicture}[x=1cm,y=1cm,every node/.style={font=\scriptsize}]

% Main parts.
\MainPart{1}{0}
\MainPart{2}{4.15}
\MainPart{3}{8.30}
\MainPart{4}{12.45}
\MainPart{5}{16.60}

% auxiliary parts.
\AuxPartNode{11}{-0.1}{-5.20}{$A_{1,1}$}
\AuxPartNode{22}{4.15}{-5.20}{$A_{2,2}$}
\AuxPartNode{41}{12.35}{-5.20}{$A_{4,1}$}
\AuxPartNode{53}{16.7}{-5.20}{$A_{5,3}$}

% Path edges between equal-colour core vertices.
\foreach \u/\v in {1/2,2/3,3/4,4/5}{
  \draw[pathEdge] (c\u a) -- (c\v a);
  \draw[pathEdge] (c\u b) -- (c\v b);
  \draw[pathEdge] (c\u c) -- (c\v c);
}

% Labels for the dummy sets 
\node[font=\tiny, text=DummyCol!88!black] at ($(D11.north)+(0,.18)$) {$D_{1,1}$};
\node[font=\tiny, text=DummyCol!88!black] at ($(D22.north)+(0,.18)$) {$D_{2,2}$};
\node[font=\tiny, text=DummyCol!88!black] at ($(D41.north)+(0,.18)$) {$D_{4,1}$};
\node[font=\tiny, text=DummyCol!88!black] at ($(D53.north)+(0,.18)$) {$D_{5,3}$};

% X_{i,t} complete to D_{i,t}
\draw[complete] (X11.north) -- (D11.south);
\draw[complete] (X22.north) -- (D22.south);
\draw[complete] (X41.north) -- (D41.south);
\draw[complete] (X53.north) -- (D53.south);

% Y_{i,t} complete to the core triple
\draw[complete] (Y11.east) .. controls +(1.18,0.38) and +(1.08,-0.38) .. (CoreBox1.east);
\draw[complete] (Y22.east) .. controls +(1.18,0.38) and +(1.08,-0.38) .. (CoreBox2.east);
%\draw[complete] (Y22.west) .. controls +(-1.18,0.38) and +(-1.08,-0.38) .. (CoreBox2.west);
\draw[complete] (Y41.east) .. controls +(1.18,0.38) and +(1.08,-0.38) .. (CoreBox4.east);
%\draw[complete] (Y53.west) .. controls +(-1.18,0.38) and +(-1.08,-0.38) .. (CoreBox5.west);
\draw[complete] (Y53.east) .. controls +(1.18,0.38) and +(1.08,-0.38) .. (CoreBox5.east);

% Z_{i,t} complete to X
\draw[complete] (X11.west) -- ++(-1.20,-0.50);
\draw[complete] (Z11.east) -- (X22.west);
\draw[complete] (Z22.east) -- (X41.west);
\draw[complete] (Z41.east) -- (X53.west);
\draw[complete] (Z53.east) -- ++(1.30,1.00);

% legend
\node[legendbox, anchor=north west] at (-1.55,-7.35) {\begin{tikzpicture}[baseline=-.5ex,every node/.style={font=\scriptsize}]
  \node[dummyset, minimum width=.26cm, minimum height=.38cm] at (0,0) {};
  \node[anchor=west] at (.24,0) {$D$};
  \node[xset, minimum width=.38cm, minimum height=.20cm] at (1.02,0) {};
  \node[anchor=west] at (1.30,0) {$X$};
  \node[yset, minimum width=.38cm, minimum height=.20cm] at (1.95,0) {};
  \node[anchor=west] at (2.23,0) {$Y$};
  \node[zset, minimum width=.38cm, minimum height=.20cm] at (2.88,0) {};
  \node[anchor=west] at (3.16,0) {$Z$};
  \draw[pathEdge] (4.05,0) -- (4.67,0);
  \node[anchor=west] at (4.82,0) {path edge};
  \draw[complete] (6.45,0) -- (7.05,0);
  \node[anchor=west] at (7.20,0) {complete-to};
\end{tikzpicture}};

\end{tikzpicture}%
}
\caption{Illustration of the construction in Theorem~\ref{thm:simulationC3Ps}. The big blue boxes depict the main parts $M_i$; they contain  $\Delta-2$ dummy sets of size $2$, as well as a triple of three vertices that represent the three possible colours of the $i$-th vertex of the path $P_s$. The green boxes depict the auxiliary parts $A_{i,j}$, split into $X$, $Y$, and $Z$ sets. The thick edges indicate a complete bipartite graph. The thin edges depict the \emph{path edges}, which can be viewed as the edges between equal colours of adjacent vertices of $P_s$.}
\label{fig:simulationC3Ps}
\end{figure}

The idea of the construction of $(G,\mathcal{P})$ in the proof of Theorem~\ref{thm:simulationC3Ps} is to ``simulate" the $3$-colour reconfiguration graph $\C_3(P_s)$ of an $s$-vertex path on $s$ vertices, by forcing it to occur as an induced subgraph of the reconfiguration graph of $(G,\mathcal{P})$.
This is inspired by~\cite[Lemma~17]{cambie2026sharp}, where a very different simulation of $C_3(P_s)$ was used for list-colouring reconfiguration. The construction is illustrated in Figure~\ref{fig:simulationC3Ps}.

\begin{proof}[Proof of Theorem~\ref{thm:simulationC3Ps}]

Put $t \coloneqq \Delta-2$ and let $s$ be a positive integer. Partition $\mathcal{P}$ will have $n \coloneqq s(\Delta-1)$  parts.

We first build $s$ main parts $M_1,\ldots,M_s$. Each main part $M_i$ consists of a core set $M_i^{\mathrm{core}} \coloneqq \{p_i^1,p_i^2,p_i^3\}$ and $t$ dummy sets $D_{i,j}$, $j\in[t]$, each of size $2$. Thus $|M_i|=3+2t=2\Delta-1$.

For each dummy set $D_{i,j}$, we introduce an auxiliary part $A_{i,j}$. This part is split into three sets
$$
A_{i,j} \coloneqq X_{i,j}\cup Y_{i,j}\cup Z_{i,j},
$$
where $|X_{i,j}|=\Delta$, $|Y_{i,j}|=1$, and $|Z_{i,j}|=\Delta-2$. Hence every auxiliary part also has size $2\Delta-1$. The partition is
$$
\mathcal P \coloneqq \{M_i:i\in[s]\}\cup\{A_{i,j}:i\in[s],\ j\in[t]\}.
$$

We now define the edges. First, between consecutive main parts, add the \emph{path edges} $p_i^c p_{i+1}^c$ for every $i\in[s-1]$ and $c\in\{1,2,3\}$. Second, make $X_{i,j}$ complete to $D_{i,j}$, for every $i,j$. Third, make $Y_{i,j}$ complete to $M_i^{\mathrm{core}}$, for every $i,j$. Finally, cyclically order all pairs $(i,j)$, and for every $(i,j)$, make $Z_{i,j}$ complete to $X_{\operatorname{succ}(i,j)}$, where $\operatorname{succ}(i,j)\neq(i,j)$ is the successor of $(i,j)$ in this cyclic order. (We note in passing that it is not essential to do this via a cyclic order; the only important thing is that each $Z_{i,j}$ is made complete to some unique $X_{p,q}$ with $(i,j)\neq (p,q)$.) There are no other edges.

We first check the maximum degree. A vertex in $M_i^{\mathrm{core}}$ has at most two path-neighbours and is adjacent to the unique vertex of $Y_{i,j}$ for each $j\in[t]$, so its degree is at most $2+t=\Delta$. A vertex in $X_{i,j}$ is adjacent to the two vertices of $D_{i,j}$ and to the $\Delta-2$ vertices of the unique $Z$-set made complete to $X_{i,j}$, so its degree is exactly $\Delta$. The unique vertex of $Y_{i,j}$ has degree $3$, and every vertex in a dummy set $D_{i,j}$ or $Z_{i,j}$ has degree exactly $\Delta$. Hence $\Delta(G)\le\Delta$.

Call an independent transversal \emph{core} if it chooses a vertex from $M_i^{\mathrm{core}}$ in every main part $M_i$, and a vertex from $X_{i,j}$ in every auxiliary part $A_{i,j}$. The crucial claim is that every core transversal can only be reconfigured to another core transversal. Indeed, in a core transversal, each vertex of $D_{i,j}$ is blocked by the selected vertex of $X_{i,j}$. The unique vertex of $Y_{i,j}$ is blocked by the selected vertex of $M_i^{\mathrm{core}}$. Finally, every vertex of $Z_{i,j}$ is blocked by the selected vertex of $X_{\operatorname{succ}(i,j)}$. Therefore the core transversals form a closed union of connected components.

Inside the core transversals, the choices in the main parts $M_1,\ldots,M_s$ are exactly proper $3$-colourings of a path on $s$ vertices: choosing $p_i^c$ corresponds to giving colour $c$ to the $i$-th path vertex, and the edges $p_i^c p_{i+1}^c$ forbid equal colours on consecutive path vertices. The auxiliary parts impose no further restrictions on these path choices, because there is no edge between the cores $\bigcup_{i} M_i^{\mathrm{core}}$ of the main parts and the cores $\bigcup_{i,j}X_{i,j}$ of the auxiliary parts.

It follows that the reconfiguration graph $\mathcal{G}$ of $(G,\mathcal{P})$ 
contains a copy $\mathcal{H}$ of the $3$-colour reconfiguration graph of the path $P_s$ as an induced 
subgraph, and moreover any shortest path between two vertices of $\mathcal{H}$ is also a shortest path in $\mathcal{G}$. In particular, the diameter of the connected component of $\mathcal{G}$ containing $\mathcal{H}$ is at least the diameter of $\mathcal{H}$. By Theorem~\ref{thm:pathreconfiguration}, $\mathcal{H}$ has diameter $\ge \frac{s^2-1}{4}$.
Plugging in $s=n/(\Delta-1)$ yields the result.
\end{proof}

\FloatBarrier

\section{Exponentially many independent transversals} \label{sec:exponentially-many}

In this section, we prove Theorem \ref{thm:exponentially-many}, demonstrating the existence of exponentially many independent transversals when parts have size $2\Delta$. We prove it in a more general form involving arbitrary part sizes $t \ge 2\Delta$ (see Theorem \ref{thm:parameterized}).

For a graph $G$ and vertex partition $\mathcal{P}$, the number of independent transversals of $(G,\mathcal{P})$ is denoted by $N(G,\mathcal{P})$.

\begin{comment}

\begin{theorem}[Wanless, Wood]
Let $G$ be a graph of maximum degree $\Delta$ and let $\mathcal{P} = \{P_1,\ldots,P_n\}$ be a partition of $V(G)$. If $|P_i| = t \ge 4\Delta$ for every $i$, then
\begin{align*}
    N(G,\mathcal{P}) \ge (t/2)^n.
\end{align*}
\end{theorem}

In fact, this was established with the \textit{maximum block average degree} $\overline{\Delta}$ in place of the maximum degree $\Delta$, the former being defined as
\begin{align*}
    \overline{\Delta} = \max_{P_i \in \mathcal{P}} \frac{1}{|P_i|} \sum_{v \in P_i} d_G(v),
\end{align*}
In this average degree formulation, the Wanless--Wood condition has been shown to be best possible for the existence of a single independent transversal. For comparison a standard application of the Lov\'asz Local Lemma gives that $N(G,\mathcal{P}) \ge \left(1 - \frac{1}{2\Delta t+1}\right)^n t^n$ whenever $|P_i| = t \ge 2e\Delta$ for every $i$. Also compare to cluster expansion Lov\'asz Local Lemma.
\end{comment}

\subsection{Double counting argument}

Before proving Theorem \ref{thm:exponentially-many}, in this subsection we give an elementary double-counting argument demonstrating exponentially many independent transversals when parts have size $2\Delta+1$. This is an adaptation of a folklore proof for showing the existence of exponentially many (list) colourings of graphs. 

\begin{prop} \label{prop:double-counting}
Let $G$ be a graph of maximum degree $\Delta$, and let $\mathcal{P}=\{P_1,\ldots,P_n\}$ be a partition of $V(G)$ with $|P_i|\ge 2\Delta$ for all $i$. Then 
\begin{align*}
    N(G,\mathcal{P}) \ge 2\Delta\prod_{i=1}^n \frac{|P_i|}{2\Delta}.
\end{align*}
In particular, if $|P_i| \ge 2\Delta + 1$ for every $i$, then 
\begin{align*}
    N(G,\mathcal{P}) \ge 2\Delta\left(1+\frac{1}{2\Delta} \right)^n.
\end{align*}
\end{prop}

\begin{proof}
For each $i$, choose a subset $Q_i\subseteq P_i$ with $|Q_i|=2\Delta$.
Write $\mathcal{Q}=\{Q_1,\ldots,Q_n\}$. The induced subgraph
$G[Q_1\cup \cdots \cup Q_n]$ has maximum degree at most $\Delta$, and its
vertex set is partitioned into $\mathcal{Q}$. By Haxell's theorem, every such choice of $\mathcal{Q}$ has at least $2\Delta$ independent transversals. Let $\mathcal{X}$ be the set of pairs $(\mathcal{Q},A)$ where $\mathcal{Q}=\{Q_1,\ldots,Q_n\}$ with $Q_i\subseteq P_i$ and $|Q_i|=2\Delta$ for every $i$, and where $A$ is an independent transversal of $\mathcal{Q}$. We will double count $|\mathcal{X}|$.

First count $\mathcal{X}$ by starting from a fixed choice of $\mathcal{Q}$. There are $\prod_{i=1}^n \binom{|P_i|}{2\Delta}$
choices for $\mathcal{Q}$, each of which has least $2\Delta$ independent transversals, so 
$|\mathcal{X}|\ge 2\Delta
\prod_{i=1}^n \binom{|P_i|}{2\Delta}.$

Now count $\mathcal{X}$ by starting from some fixed choice of independent transversal $A$ of $\mathcal{P}$. For $A$ to be an independent transversal of
$\{Q_1,\ldots,Q_n\}$, each $Q_i$ must contain the unique vertex of
$A\cap P_i$. After this vertex has been included, the remaining
$2\Delta-1$ vertices of $Q_i$ may be chosen arbitrarily from the other
$|P_i|-1$ vertices of $P_i$. Hence, the number of choices of $\mathcal{Q}$ compatible with $A$ is
$\prod_{i=1}^n \binom{|P_i|-1}{2\Delta-1}.$
Thus
$|\mathcal{X}|=
N\prod_{i=1}^n \binom{|P_i|-1}{2\Delta-1},$ where $N = N(G,\mathcal{P})$.

Combining the two estimates yields
$N\prod_{i=1}^n \binom{|P_i|-1}{2\Delta-1}
= |\mathcal{X}| \ge 2\Delta
\prod_{i=1}^n \binom{|P_i|}{2\Delta}$.
Therefore, 
\begin{align*}
    N(G,\mathcal{P}) = N\ge 2\Delta\prod_{i=1}^n\frac{\binom{|P_i|}{2\Delta}}{\binom{|P_i|-1}{2\Delta-1}}= 2\Delta\prod_{i=1}^n \frac{|P_i|}{2\Delta}
\end{align*}
as required.
\end{proof}

The above proposition and proof, though not fully optimized, easily adapt to other settings. For example, if $G$ is a chordal graph, then a result of Aharoni, Berger, and Ziv \cite{aharoni2002tree} guarantees that every vertex lies in an independent transversal whenever $|P_i| \ge \Delta+1$ for every $i$. The above proof then yields that, when $G$ is chordal and $|P_i| \ge \Delta+2$ for every $i$, we have $N(G,\mathcal{P}) \ge (\Delta+1)\left(1+\frac{1}{\Delta+1}\right)^n$.

We also mention that when $n$ is relatively small ($2\Delta$ or smaller), then we obtain a better lower bound on $N(G,\mathcal{P})$ by applying the above proof on another refined version of Haxell's Theorem \ref{thm:existence}, due to Haxell and Szab\'o \cite{haxell2006odd}, stated below.

\begin{theorem}[Haxell, Szab\'o \cite{haxell2006odd}]
Let $G$ be a graph of maximum degree at most $\Delta$ with a vertex partition $\mathcal{P}=\{P_1,\ldots, P_n\}$. Then $(G,\mathcal{P})$ has an independent transversal whenever
\begin{itemize}
    \item $|P_i| \ge \left\lfloor\frac{2\Delta (n-1)}{n} \right\rfloor+1$ for every $i$, if $n \ge 2$ is even;
    \item $|P_i| \ge \left\lfloor\frac{2\Delta (n-2)}{n-1} \right\rfloor+1$ for every $i$, if $n \ge 3$ is odd.
\end{itemize}
\end{theorem}

However, the resulting lower bound is superseded by the following theorem, mentioned in the Introduction, of Tang and Zhao \cite{tang2025number} when $n \ge 4$ is even and $\Delta$ is sufficiently large.

\begin{theorem}[Tang, Zhao \cite{tang2025number}] \label{thm:tang-zhao}
For any even integer $n \ge 4$, there exist $\Delta_0 \ge 1$ and $\lambda_0 > 0$ such that whenever $\Delta \ge \Delta_0$ and $0 \le k \le \lambda_0 t$ and $t \ge \left\lfloor\frac{2(\Delta+k)(n-1)}{n} \right\rfloor+1$, the following holds. If $G$ is a graph of maximum degree $\Delta$, $\mathcal{P} = \{P_1, \ldots, P_n\}$ is a partition of $V(G)$, and $|P_i| = t$ for every $i$, then
\begin{align*}
    N(G,\mathcal{P}) \ge \frac{1}{4}(2k(n-1) + 1)\left(\frac{t n}{2(n-1)} \right)^{n-1}.
\end{align*}
In particular:
\begin{itemize}
    \item For any even $n \ge 4$ and sufficiently large $\Delta$, if $|P_i| \ge \left\lfloor \frac{2\Delta(n-1)}{n} \right\rfloor+1$ for every $i$, then $N(G,\mathcal{P}) \ge \frac{1}{4}\Delta^{n-1}$.
    \item For any $n \ge 3$ and sufficiently large $\Delta$, if $|P_i| \ge 2\Delta$ for every $i$, then $N(G,\mathcal{P}) \ge c_n\Delta^{n}$ for some constant $c_n > 0$ depending only on $n$.
\end{itemize}
\end{theorem}

\subsection{Isoperimetric approach}

In this subsection, we prove the following more general version of Theorem \ref{thm:exponentially-many}, as an application of edge-isoperimetric inequalities. 

\begin{theorem} \label{thm:parameterized}
Let $G$ be a graph of maximum degree $\Delta$, and let $\mathcal{P}=\{P_1,\ldots,P_n\}$ be a partition of $V(G)$.
For every integer $t$ satisfying
$2\Delta\le t\le\min_i |P_i|$, we have
$$
N(G,\mathcal P)\ge
\frac{\prod_{i=1}^n |P_i|}{t^{n\Delta/(t-1)}}. % \wouter{\frac{\prod_{i=1}^n |P_i|}{t^{n\Delta/(t-1)}}}
$$
In particular: %(\wouter{check consequences for the two corollaries})
\begin{itemize}
    \item If $|P_i| \ge 2\Delta+1$ for every $i$, then $N(G,\mathcal{P}) \ge (2\Delta+1)^{\frac{n}{2}}$.
    \item If $|P_i| \ge 2\Delta$ for every $i$, then $N(G,\mathcal{P}) \ge (2\Delta)^{\left(1-\frac{1}{2\Delta-1}\right)\frac{n}{2}}$.
\end{itemize}
\end{theorem}

The basic idea is to view our reconfiguration graph on independent transversals as a subgraph of the Hamming graph, show that it induces many edges, and thus by isoperimetric inequalities deduce that it has many vertices. 

For integers $n,t \ge 1$, the \textit{$t$-ary $n$-dimensional Hamming graph} $H(n,t)$ has vertex set $[t]^n$, with two words in $[t]^n$ connected by an edge if they differ in exactly one coordinate. We use the following formulation of Lindsey's edge-isoperimetric inequality on the Hamming graph \cite{lindsey1964assignment}
(see also inequality (3) in \cite{diskin2025isoperimetry}, from which it easily follows).
%Much later, Squier, Torence, Vogt  independently proved it explicitly as their proposition 2 in https://doi.org/10.1016/S0893-9659(01)80030-6 }
\begin{lem}
\label{lem:hamming-entropy}
For any vertex subset $\mathcal{S}$ of the Hamming graph $H(n,t)$, with $n\ge 1$ and $t \ge 2$, the number of edges induced by $\mathcal S$ is at most
\begin{equation*}
    \frac{1}{2}(t-1)|\mathcal{S}|\log_t|\mathcal S|.
\end{equation*}
\end{lem}

\begin{proof}[Proof of Theorem \ref{thm:parameterized}]
We prove the stated lower bound when every part has size exactly $t$, where the claimed lower bound is $t^{n} \cdot t^{-n\Delta/(t-1)}$. To get the statement for general part sizes $|P_i|$, we then apply the same double counting argument as in Proposition~\ref{prop:double-counting}.

Identify each part $P_i$ with the set $[t]$ in arbitrary fashion, and let $\mathcal T \subseteq [t]^n$ represent the family of independent transversals of $(G,\mathcal{P})$ under these identifications. We view $\mathcal{T}$ as a vertex subset of the Hamming graph $H(n,t)$. Since $t\ge2\Delta$, Haxell's Theorem \ref{thm:existence} gives $\mathcal T\ne\varnothing$.

Fix an independent transversal $A=\{a_1,\ldots,a_n\}\in\mathcal T$, with $a_i\in P_i$ for every $i$, identified as vertices of $H(n,t)$. Define the set $P_i(A) = \{x\in P_i:(A\setminus\{a_i\})\cup\{x\}\in\mathcal T\}$. Every vertex in $P_i \setminus P_i(A)$ has a neighbour in $A \setminus \{a_i\}$, which implies that $\sum_{i=1}^{n} (|P_i|-|P_i(A)|) \le \sum_{j=1}^{n} d_G(a_j) \le n\Delta$.
%so \textcolor{red}{$|P_i| - |P_i(A)| \le \sum_{j \neq i} d_G(a_j) \le (n-1)\Delta$}.
%\wouter{I believe that instead we need an estimate on the \emph{sum} over all $i$. More precisely, for the same reason we have $\sum_{i=1}^{n} |P_i|-|P_i(A)| \le \sum_{j=1}^{n} d_G(a_j) \le n\Delta.$ As a consequence, all red-colored $+\Delta$s below then need to be removed.}
Thus, the degree of $A$ in the subgraph of $H(n,t)$ induced by $\mathcal{T}$ satisfies
%\begin{align*}
%    d_{\mathcal{T}}(A) &= \sum_{i=1}^n (|P_i(A)|-1) = n(t-1) - \sum_{i=1}^n (|P_i| - |P_i(A)|) \\
%    &\ge n(t-1) - (n\textcolor{red}{-1})\Delta \ge n(t - 1 - \Delta)\textcolor{red}{+\Delta}.
%\end{align*}
\begin{align*}
    d_{\mathcal{T}}(A) &= \sum_{i=1}^n (|P_i(A)|-1) = n(t-1) - \sum_{i=1}^n (|P_i| - |P_i(A)|) \\
    &\ge n(t-1) - n\Delta \ge n(t - 1 - \Delta).
\end{align*}
Since $A$ was an arbitrary vertex in $\mathcal{T}$, we get that the number of edges induced by $\mathcal{T}$ in $H(n,t)$ is at least $\frac{1}{2}n(t - 1 - \Delta)|\mathcal{T}|$. Combining this with Lemma~\ref{lem:hamming-entropy} gives $\frac{1}{2}n(t - 1 - \Delta) |\mathcal{T}| \le \frac{1}{2}(t-1)|\mathcal{T}|\log_t|\mathcal T|$, and therefore
\begin{equation*}
|\mathcal T|\ge t^{n(t-1-\Delta)/(t-1)} = t^n \cdot t^{-n\Delta/(t-1)},
\end{equation*}
as required. 
%\wouter{The final equality then becomes $t^n \cdot t^{-n(\Delta/(t-1))}$}
\end{proof}

\section{Conjecture~\ref{conj:supersaturation} for small number of parts} \label{sec:small-number-of-parts}

In this section, we prove the first few cases of Conjecture~\ref{conj:supersaturation}, namely when $n \le 3$. Recall that $N(G,\mathcal{P})$ denotes the number of independent transversals of $(G,\mathcal{P})$, and recall that Conjecture \ref{conj:supersaturation} of Buys, Kang, and Ozeki \cite{buys2025reconfiguration} states that, if $|P_i| = 2\Delta$ for all $i$, then
\begin{align*}
    N(G,\mathcal{P}) \ge 2\Delta^n.
\end{align*}

This lower bound would be best possible, as it is attained by the following basic cyclic construction, a special case of Wdowinski's characterization \cite{wdowinski2026tight} of when the reconfiguration graph of $(G,\mathcal{P})$ is disconnected. 

\begin{example} \label{ex:cyclic}
Let $K^1, \ldots, K^n$ be disjoint copies of $K_{\Delta,\Delta}$, where $K^i$ has bipartition $A_i \cup B_i$. Let $G = K^1 \cup \cdots \cup K^n$, and let $\mathcal{P} = \{P_1, \ldots, P_n\}$ consist of $P_i = A_{i-1} \cup B_{i}$ for every $i \in \{1, \ldots, n\}$ (where $A_0 = A_n$). Then every independent transversal of $(G, \mathcal{P})$ is either a transversal of $A_1, \ldots, A_n$ or a transversal of $B_1, \ldots, B_n$, and there are $2 \Delta^n$ such transversals. 
\end{example}

Observe that Conjecture~\ref{conj:supersaturation} trivially holds when $n = 1$, with equality. Next, it is straightforward to show that the conjecture holds when $n=2$, although in this case the cyclic construction of Example \ref{ex:cyclic} is not the unique minimizer.

\begin{prop} \label{prop:two-parts}
Let $G$ be a graph of maximum degree $\Delta$ and let vertex partition $\mathcal{P} = \{P_1, P_2\}$ satisfy $|P_1|, |P_2| \ge 2\Delta$. Then $N(G,\mathcal{P}) \ge 2\Delta^2$, and equality holds if and only if $|P_1| = |P_2| = 2\Delta$ and $G$ is a $\Delta$-regular bipartite graph with bipartition $\{P_1, P_2\}$.
\end{prop}

\begin{proof}
The independent transversals of $(G,\mathcal{P})$ consist of a vertex $v \in P_1$ together with a non-neighbour of $v$ in $P_2$. Every vertex in $P_1$ has at least $|P_2| - \Delta$ non-neighbours in $P_2$, so all together there are at least $|P_1|(|P_2| - \Delta) \ge 2\Delta \cdot \Delta = 2\Delta^2$ independent transversals of $(G, \mathcal{P})$. Equality holds if and only if $|P_1| = |P_2| = 2\Delta$ and every vertex of $P_1$ has exactly $\Delta$ neighbours in $P_2$. In this case, $G$ is a $\Delta$-regular bipartite graph with bipartition $\mathcal{P}$.
\end{proof}

Our next goal is to prove the conjecture for $n = 3$, and demonstrate that the described cyclic construction of Example \ref{ex:cyclic} is indeed the unique minimizer of $N(G,\mathcal{P})$. 

\begin{theorem} \label{thm:three-parts}
Let $G$ be a graph of maximum degree $\Delta$ and let vertex partition $\mathcal{P} = \{P_1, P_2, P_3\}$ satisfy $|P_1|, |P_2|, |P_3| \ge 2\Delta$. Then $N(G,\mathcal{P}) \ge 2\Delta^3$, and equality holds if and only if $(G, \mathcal{P})$ is a cyclic construction of the form described in Example \ref{ex:cyclic}.
\end{theorem}

We prove Theorem \ref{thm:three-parts} using the following lemma, which comes from an inclusion-exclusion argument. Recall that a \textit{transversal} of partition $\mathcal{P}$ is a set with one element from each part of $\mathcal{P}$.

\begin{lem} \label{lem:inclusion-exlcusion}
For a graph $G$ and partition $\mathcal{P}=\{P_1,\ldots, P_n\}$ of $V(G)$ with $|P_i|= t$ for all $i$, we have
\begin{align*}
    N(G,\mathcal{P}) \ge t^{n-1} \left( t - \frac{\Delta n}{2} \right)  + \sum_{\substack{T~ \mathrm{transversal} \\ \mathrm{of}~\mathcal{P}}} \max\{0, |E(T)|-1\}.
\end{align*}
\end{lem}

\begin{proof}
For an arbitrary transversal $T=\{v_1,\ldots,v_n\}$, with $v_i \in P_i$ for all $i$, let $X_{i,j}$ be the indicator variable for $v_i, v_j$ inducing an edge. Then $T$ is an independent transversal of $(G,\mathcal{P})$ precisely when
\begin{equation*}
    \prod_{\{i,j\} \subseteq [n]} (1 - X_{i,j})
\end{equation*}
equals $1$, and otherwise $X_T$ equals $0$. Thus, the number of independent transversals of $(G,\mathcal{P})$ is 
\begin{equation} \label{eq:number-product}
    N(G,\mathcal{P}) = \sum_{\substack{T = \{v_1,\ldots,v_n\} \\ \text{transversal}}} \prod_{\{i,j\} \subseteq [n]} (1 - X_{i,j}).
\end{equation}

%We may assume that $|P_i| = t$ for all $i$ (since deleting vertices does not increase the number of independent transversals) and that the $P_i$'s are independent sets of $G$ (since deleting edges within the $P_i$'s does not change the number of independent transversals).
%
%this reduction was not quite true, due to the correction term in the lower bound. To correct for that, the lemma has now only been stated for parts of size precisely $t$, which is all we need in the application.
For a fixed transversal $T = \{v_1,\ldots,v_n\}$, we expand the associated product as
\begin{equation*}
    \prod_{\{i,j\} \subseteq [n]} (1 - X_{i,j}) = 1 - \sum_{\{i,j\} \subseteq [n]} X_{i,j} + \sum_{\substack{\{i_1,j_1\}, \{i_2,j_2\} \subseteq [n] \\ \text{distinct}}} X_{i_1,j_1} X_{i_2,j_2} - \cdots
\end{equation*}
Summing over all $T$, the first term becomes $t^n$, the number of transversals $T$ of $\mathcal{P}$. The second term becomes $-t^{n-2}|E(G)|$ because a given edge lies in $t^{n-2}$ transversals $T$. Equation (\ref{eq:number-product}) becomes
\begin{align*}
    N(G,\mathcal{P}) &= t^n - t^{n-2}|E(G)| + \sum_{T \text{ transversal}} \sum_{k=2}^{\binom{n}{2}} (-1)^k \binom{|E(T)|}{k} \nonumber \\
    &= t^n - t^{n-2}|E(G)|  + \sum_{T \text{ transversal}} \max\{0, |E(T)|-1\} \nonumber \\
    &\ge t^{n-1} \left( t - \frac{\Delta n}{2} \right)  + \sum_{T \text{ transversal}} \max\{0, |E(T)|-1\},
\end{align*}
where the last inequality uses that $|E(G)| \le \frac{\Delta t n}{2}$. This gives the required result.
\end{proof}

\begin{proof}[Proof of Theorem \ref{thm:three-parts}]
The lower bound $N(G,\mathcal{P}) \ge 2\Delta^3$ comes directly from Lemma \ref{lem:inclusion-exlcusion} by plugging in $t = 2\Delta$ and $n = 3$ (the last sum is always nonnegative). Equality holds if and only if $G$ is $\Delta$-regular and $|E(T)| \in \{0,1\}$ for every transversal $T$ of $\mathcal{P}$, and each $P_i$ is an independent set. In this case, we have that the non-neighbours of $P_1$ inside $P_2 \cup P_3$ form a $\Delta$-regular bipartite subgraph of $G$ with bipartition (inherited from) $P_2 \cup P_3$. Since $|N(P_1)| \ge 2\Delta$, there are at most $2\Delta$ non-neighbours of $P_1$ inside $P_2 \cup P_3$. Hence by $\Delta$-regularity these non-neighbours must induce a $K_{\Delta,\Delta}$. This then implies that $|N(P_1)| = 2\Delta$. Moreover, $N(P_1)$ is an independent set, and no vertex in $P_1$ can have neighbours in both $P_2$ and $P_3$ (by $|E(T)| \in \{0,1\}$ for every transversal $T$). From these, we obtain that $P_1 \cup (N(P_1) \cap P_2)$ and $P_1 \cup (N(P_1) \cap P_3)$ contain vertex-disjoint induced copies of $K_{\Delta,\Delta}$, with those bipartitions. This gives the cyclic construction $(G,\mathcal{P})$.
\end{proof}

\section*{Declaration of use of AI}
The majority of this work was done in May, June and July 2026.
For Theorems \ref{thm:parameterized} and \ref{thm:three-parts}, the main ideas were originally found through multiple queries of ChatGPT-5.5.  The authors verified, cleaned up, and generalized the details. ChatGPT-5.5 was also used for initial versions of Figure~\ref{fig:simulationC3Ps}.

\section*{Acknowledgments}

We would like to thank Pjotr Buys, Carla Groenland, and Jan van den Heuvel for fruitful discussions. This project started during the ``Phase Transitions in Combinatorics, Algorithms, Probability'' (PhaseCAP) Research Semester Programme at CWI in 2026. W.CvB. was supported by the Belgian National Fund for Scientific Research (FNRS). P.H. was partially supported by NSERC as well as the CWI Visiting Programme. R.J.K. was partially supported by the Dutch Research Council (NWO) grants OCENW.M20.009 and OCENW.M.25.085. R.W. was supported in part by the Austrian Science Fund (FWF) [10.55776/F1002]. For open access purposes, the authors have applied a CC BY public copyright license to any author accepted manuscript version arising from this submission.

\bibliographystyle{habbrv}
\bibliography{reconfiguration}

\end{document}